\documentclass[12pt]{article}
\usepackage{stmaryrd}
\usepackage{mathrsfs}
\usepackage{amsmath}
\usepackage{amsmath,amsthm,amssymb}
\usepackage{amsfonts}
\usepackage{latexsym}
\date{}

\begin{document}
\title{Partial tilting modules over $m$-replicated algebras}
\author{{\small Shunhua Zhang}\\
{\small  Department of Mathematics,\ Shandong University,\ Jinan
250100, P. R. China }\\
{\small   {\it Email addresses}: shzhang@sdu.edu.cn   }}

\pagenumbering{arabic}

\maketitle
\begin{center}
 \begin{minipage}{120mm}
   \small\rm
   {\bf  Abstract}\ \ Let $A$ be a hereditary algebra over an algebraically
   closed field $k$ and $A^{(m)}$ be the $m$-replicated algebra of
   $A$. Given an $A^{(m)}$-module $T$, we denote by $\delta (T)$ the number of
   non isomorphic indecomposable summands of $T$.  In this paper,  we prove that
   a partial tilting $A^{(m)}$-module $T$ is a tilting $A^{(m)}$-module if and only if
   $\delta (T)=\delta (A^{(m)})$, and that every partial tilting $A^{(m)}$-module has complements.
   As an application, we deduce that the tilting quiver $\mathscr{K}_{A^{(m)}}$ of $A^{(m)}$
   is connected. Moreover, we investigate the number
   of complements to almost tilting modules over
   duplicated algebras.

   \vskip 0.1in


\end{minipage}
\end{center}

\section {Introduction}

\vskip 0.2in

  Let $A$ be  an  Artin algebra. We denote by $A$-mod the category of all finitely
generated left $A$-modules, and we always assumed that subcategories
of $A$-modules closed under isomorphisms and direct summands. We
denote by ${\rm pd}_A X$ the projective dimension of an $A$-module
$X$ and by ${\rm gl.dim}\ A$ the global dimension of $A$. Given an
$A$-module $M$, we denote by ${\rm add} \ M$ the full subcategory
having as objects the direct sums of indecomposable summands of $M$
and by $\delta (M)$ the number of non isomorphic indecomposable
summands of $T$.

\vskip 0.2in

  A module $T\in A$-mod is called a (generalized) tilting module if
the following conditions are satisfied:\\
{(1)} ${\rm pd}_A T =n<\infty$;\\
{(2)} ${\rm Ext}_A^{i}(T,T)=0 $ for all $i >0$;\\
{(3)} There is  a long exact sequence
 $$0\longrightarrow
A\longrightarrow T_{0}\longrightarrow T_{1}
\longrightarrow\cdot\cdot\cdot\longrightarrow T_{n}\longrightarrow
0$$ with $T_{i}\in {\rm add}\ T $ for $0\leq i\leq n$.

\vskip 0.2in

  An $A$-module $M$  satisfying the conditions $(1)$ and $(2)$ of
the definition above is called a partial tilting module and if
moreover $\delta(M)=\delta(A)-1 $, then  $M$ is called an almost
complete tilting module. Let $M$ be a partial tilting module and $X$
be an $A$-module such that $M\oplus X$ is a tilting module and ${\rm
add}M\cap {\rm add}X=0$. Then X will be called a complement to $M$.

\vskip 0.2in

It is well known that in the classical situation $M$ always admits a
complement and $M$ is a tilting module if and only if $\delta (M) =
\delta (A)$ [5]. However, in general situations complements do not
always exist, as shown in [13]. Moreover it is an important open
problem whether for a partial tilting module $M$ with $\delta (M) =
\delta (A)$ is sufficient for $M$ to be a tilting module. In this
paper, we prove that this is true for $m$-replicated algebra
$A^{(m)}$, i.e., a partial tilting $A^{(m)}$-module $T$  is  tilting
if and only if $\delta (T)=\delta (A^{(m)})$, and  every partial
tilting $A^{(m)}$-module has complements.

\vskip 0.2in

Let $A$ be a hereditary algebra over an algebraically closed field
$k$ and $A^{(m)}$ be the $m$-replicated algebra $A$. This kind of
algebras, introduced in [1,2], gives a one-to-one correspondence
between basic tilting  $A^{(m)}$-modules with projective dimension
at most $m$ and basic  tilting objects in $m-$cluster category
$\mathscr{C}_m(A)$ and it is proved that a faithful partial tilting
$A^{(m)}$-module $T$ with pd$_{A^{(m)}}T\leq m$ is tilting if
$\delta (T)=\delta (A^{(m)})$. In [12], we proved that the
presentation dimension of $A^{(m)}$ is at most 3, and in [11], we
investigated complements to the almost complete tilting
$A^{(m)}$-modules and proved that a faithful almost complete tilting
$A^{(m)}$-module $T$ with pd$_{A^{(m)}}T\leq m$ has  exactly $m+1$
indecomposable non-isomorphic complements of projective dimensions
at most $m$. This motivates the further investigation on the
generalized partial tilting modules over $m$-replicated algebras.

\vskip 0.2in

 Now we state our main results of this paper in the following theorems.

\vskip 0.2in

{\bf Theorem 1.} {\it A partial tilting $A^{(m)}$-module $T$ is
tilting if and only if $\delta (T)= \delta (A^{(m)})$.}

\vskip 0.2in

{\bf Theorem 2.} {\it  Let $M$ be a partial tilting
$A^{(m)}$-module. Then $M$ admits a complement $C$ in
$A^{(m)}$-mod.}

\vskip 0.2in

Note that Theorem 1 and Theorem 2 above generalized the results in
[2], here the restriction on projective dimension for partial
tilting modules over $m$-replicated algebra $A^{(m)}$ is removed.
Moreover, our proofs are very deferent from [2].

\vskip 0.2in

 {\bf Theorem 3.} {\it The tilting quiver
 $\overrightarrow{\mathscr{K}}_{A^{(m)}}$ of
 $A^{(m)}$ is connected.}

\vskip 0.2in

This paper is arranged as the following. In section 2, we fix the
notations and recall some necessary facts needed for our further
research. Section 3 is devoted to the proof of Theorem 1, Theorem 2
and Theorem 3. In section 4, we investigate the number of
complements to an almost tilting module over duplicated algebras.

\vskip 0.2in

\section {Preliminaries}

\vskip 0.2in

Let $A$ be a finite dimensional algebra over an algebraically closed
field $k$. We denote by $A$-mod the category of all finitely
generated left $A$-modules and by $A$-ind  the full subcategory of
$A$-mod containing exactly one representative of each isomorphism
class of indecomposable $A$-modules. $D={\rm Hom}_k(-,\ k)$ is the
standard duality between $A$-mod and $A^{op}$-mod, and $\tau_A$ is
the Auslander-Reiten translation of $A$. The Auslander-Reiten quiver
of $A$ is denoted by $\Gamma_A$.

\vskip 0.2in

Let $\mathcal{C}$ be a full subcategory of $A$-mod,
$C_{M}\in\mathcal{C}$ and $\varphi :C_M\longrightarrow M$ with
$M\in$ $A$-mod. The morphism $\varphi$ is a right
$\mathcal{C}$-approximation of $M$ if the induced  morphism ${\rm
Hom}_A(C,C_{M})\longrightarrow {\rm Hom}_A(C,M)$ is surjective for
any $C\in\mathcal{C}$. A minimal right $\mathcal{C}$-approximation
of $M$ is a right $\mathcal{C}$-approximation which is also a right
minimal morphism, i.e., its restriction to any nonzero summand is
nonzero. The subcategory $\mathcal{C}$ is called contravariantly
finite if any module $M\in$ $A$-mod admits a (minimal) right
$\mathcal{C}$-approximation. The notions of (minimal) left
$\mathcal{C}$-approximation and of covariantly finite subcategory
are dually defined. It is well known that add $M$ is both a
contravariantly finite subcategory and a covariantly finite
subcategory.

\vskip 0.2in

Given any module $M\in A $-mod, we denote by $M^{\bot}$ the
subcategory of $A$-mod with objects $X\in A$-mod satisfying ${\rm
Ext}_A^{i}(M,X)=0$ for all $i\geq 1 $  and by $^{\perp}M$ the
subcategory of $A$-mod with objects $X\in A$-mod satisfying ${\rm
Ext}_A^{i}(X,M)=0$ for $ i\geq 1$. We denote by $\Omega_{A}^{i}M$
 and $\Omega_{A}^{-i}M$ the i-th syzygy and cosyzygy  of $M$ respectively,
 and denote by gen $M$
the subcategory of $A$-mod whose objects are generated by $M$. We
may decompose $M$ as $M\cong\oplus_{i=1}^{m}M_{i}^{d_{i}}$, where
each $M_{i}$ is indecomposable, $d_{i}>0$ for each $i$, and $M_{i}$
is not isomorphic to $M_{j}$ if $i\neq j$. The module $M$ is called
basic if $d_{i}=1 $ for any $i$. The number of non-isomorphic
indecomposable modules occurring in the direct sum decomposition
above is uniquely determined and it is denoted by $\delta(M)$.

\vskip 0.2in

Let $M, N $ be two indecomposable $A$-modules. A path from $M$ to
$N$ in $A$-ind is a sequence of non-zero morphisms
  $$M=M_0\stackrel{f_1} \longrightarrow M_1\stackrel{f_2} \longrightarrow\cdots
  \stackrel{f_t} \longrightarrow M_t =N$$
  with all $M_i$ in  $A$-ind. Following [14], we denote the
  existence of such a path by $M\leq N$. We say that $M$ is a
  predecessor of $N$ (or that $N$ is a successor of $M$).

\vskip 0.2in

 From now on, let $A$ be a finite dimensional hereditary
 algebra over algebraically closed field $k$. According to [2], we
  define the m-replicated algebra of $A$ as the (finite dimensional)
  matrix algebra

$$
A^{(m)}=\left [\begin{array}{cccccc}
   A_0&0 & \cdots& \cdots &\cdots & 0  \\
    Q_1& A_1 &0 &\cdots& \cdots&  0  \\
  0& Q_2&A_2&0 &\cdots & 0 \\
   \vdots & \ddots &\ddots \\
   0& \cdots& & 0&Q_m&A_m \end{array}\right ] \\
$$
where $A_i=A,\ Q_i=DA$ for all i and all the remaining coefficients
are zero and multiplication is induced from the canonical
isomorphisms $A\otimes_A DA\cong\! _{A}DA_{A}\cong DA\otimes_A A $
and the zero morphism $DA\otimes_A DA\longrightarrow 0$.

\vskip 0.2in

We identify $A$ with $A_0$ and each $A_i$-${\rm ind}$ with the
corresponding full convex subcategory of $A^{(m)}$-${\rm ind}$
 for $0\leq i\leq m$.  We denote by $\Sigma_0$ the set of all non-isomorphic
indecomposable
 projective $A_0$-modules and denote $\Omega_{A^{(m)}}^{-i}\,\Sigma_0$ by $\Sigma_i$.

\vskip 0.2in

\vskip 0.1in If $m=1$, then $A^{(1)}$ is the duplicated algebra of
$A$ (see [1]). Also from [2], we have that $m+1 \leq {\rm gl.dim }\
A^{(m)}\leq 2m+1$. Moreover, if $A$ is representation-infinite, then
${\rm gl.dim }\ A^{(m)}= 2m+1$.

\vskip 0.2in

Let $A'$ be the right repetitive algebra of $A$ defined in [1,2].
The next lemma also is proved  in [2].

\vskip 0.1in

{\bf Lemma 2.1.}\ {\it {\rm (1)}The standard embeddings $A_{i}-{\rm
ind} \hookrightarrow A^{(m)}-{\rm ind }$ (where $0\leq i\leq m$) and
$A^{(m)}-{\rm ind } \hookrightarrow A'-{\rm ind }$ are full, exact,
preserve indecomposable modules, almost split sequences and
irreducible morphisms.

\vskip 0.1in

{\rm(2)}\ Let $M$ be an indecomposable $A'$-module which is not
projective and $k\geq 1$. Then {\rm pd} $M =k$ if and only if
 $\Sigma_{k-1}< M\leq \Sigma_{k}.$

\vskip 0.1in

 {\rm(3)}\ Let $M$ be an indecomposable $A^{(m)}$-module
which is not in {\rm ind} $A^{(0)}$. Then its projective cover in
$A^{(m)}$-{\rm mod}  is projective-injective and coincides with its
projective cover in $A'$-{\rm mod}.}

\vskip 0.2in

The following lemma is the main results in [11].

 \vskip 0.2in

{\bf Lemma 2.2.}  {\it Let $T$ be a faithful almost complete tilting
$A^{(m)}$-module with  ${\rm pd} _{A^{(m)}}T\leq m$. Then there
exists an exact sequence
  $$ (*) \ \ \ \  0\longrightarrow
X_0\stackrel{g_{0}}\longrightarrow
T_{0}^\prime\stackrel{g_{1}}\longrightarrow T_{1}^\prime
\longrightarrow\cdot\cdot\cdot\longrightarrow
T_{m-2}^\prime\stackrel{g_{m-1}}\longrightarrow
T_{m-1}^\prime\longrightarrow X_m\longrightarrow 0$$ in
$A^{(m)}$-mod,  such that

\vskip 0.1in

{(1)} $T_i'\in{\rm add}\, T $ for all $0\leq i\leq m-1$,

\vskip 0.1in

{(2)} $ X_i={\rm Coker}\ g_{i-1}$ for $1\leq i\leq m$, $i\leq {\rm
pd}_{A^{(m)}}X_i\leq i+1$
 for  $0\leq i\leq m$,

\vskip 0.1in

{(3)} each of the induced monomorphisms $X_i\hookrightarrow
T_i^\prime$ is a  minimal left  ${\rm add}\ T$-approximation.

\vskip 0.1in

 {(4)}  $T$ has exactly $m+1$ non-isomorphic indecomposable
 complements $X_0, X_1,$ $\cdots$, $ X_m$ with  projective dimensions at most $m$.

 \vskip 0.1in

 {(5)} If ${\rm pd} _{A^{(m)}}T =t \leq m$, the complements to
$T$ have the distribution as the following.

 ($i$) If ${\rm pd}_{A^{(m)}} X_0=0$, then ${\rm pd}_{A^{(m)}}X_i=i$ for each i.

($ii$) If ${\rm pd}_{A^{(m)}}  X_0\neq 0$, then there exists a
unique j with $0\leq j\leq t-1 $ such that ${\rm pd}_{A^{(m)}}X_j =
{\rm pd}_{A^{(m)}}X_{j+1}=j+1$, ${\rm pd}_{A^{(m)}}X_i=i+1$ for
$0\leq i\leq j$ and ${\rm pd}_{A^{(m)}} X_i=i$ for $i\geq j+1$.

 \vskip 0.1 in

 {(6)} The number of non-isomorphic indecomposable complements
to $T$ is either $2m+1$ or $2m+2$.}

 \vskip 0.2in

According to Proposition 2.2 in [6], we have the following Lemma.

 \vskip 0.2in

{\bf Lemma 2.3.}  {\it Let $M$ be a partial tilting
$A^{(m)}$-module. If $^{\bot}M$ is contravariantly finite
subcategory in $A$-mod, then $M$ admits a complement $C$.}

 \vskip 0.2 in

 The following Lemma is proved in [9].

\vskip 0.2 in

 {\bf Lemma 2.4.}  {\it Let $M$ be an almost complete
tilting module with an indecomposable complement $X$.

{(1)} If $X$ is generated by $M$ and $f:M'\longrightarrow X$ is a
minimal right ${\rm add}\, M$-approximation of $X$, then ${\rm Ker}
f$ is also an indecomposable complement not isomorphic to $X$ and
${\rm Ker}f\longrightarrow M'$ is a minimal left ${\rm add}\,
M$-approximation of ${\rm Ker} f$.

{(2)} If $X$ is cogenerated by $M$ and $g:X\longrightarrow M''$ is a
minimal left ${\rm add}\, M$-approximation of $X$, then ${\rm
Coker}\ g$ is also an indecomposable complement not isomorphic to
$X$ and $M''\longrightarrow {\rm Coker}\ g$ is a minimal right ${\rm
add}\, M$-approximation of ${\rm Coker}\ g$.}

\vskip 0.2in

For the Auslander-Reiten quivers of $A'$ and $A^{(m)}$, we refer to
[2]. Throughout this paper, we follow the standard terminology and
notation used in the representation theory of algebras, see [4, 14].

\vskip 0.2in

\section {Partial tilting $A^{(m)}$-modules}

\vskip 0.2in

In this section, we prove Theorem 1, Theorem 2 and Theorem 3
promised in the introduction.

\vskip 0.2in

{\bf Theorem 3.1.} {\it Let $T$ be a partial tilting
$A^{(m)}$-module.  Then $T$ is a tilting $A^{(m)}$-module if and
only if $\delta (T)= \delta (A^{(m)})$.}

\vskip 0.1in

{\bf Proof.}  We only need to prove that $T$ is a tilting
$A^{(m)}$-module whenever $\delta (T)= \delta (A^{(m)})= n(m+1)$.
Now we suppose ${\rm pd}\ T = t$.

\vskip 0.1in

Case I. Let $t \leq m$. By [2], we know that $T$ admits a complement
$C$ in $A^{(m)}$-mod, i.e., $T\oplus C$ is a tilting
$A^{(m)}$-module and $\delta (T)= \delta (A^{(m)})=(m+1)\delta
(A)=\delta (T\oplus C)$. We have that ${\rm add} \ C\subset {\rm
add} \ T$, hence $T=T\oplus C$ is a tilting $A^{(m)}$-module.

\vskip 0.1in

Case II. Let $t > m$.   Then $T$ can be regarded as a partial
tilting $A^{(t)}$-module.  As in Case I and by [2], we know that $T$
admits a complement $X$ in $A^{(t)}$-mod, i.e., $T\oplus X$ is a
tilting $A^{(t)}$-module and $\delta (T\oplus X)=\delta (A^{(t)})=
(t+1)\delta (A)$.

 We have that $\delta (X)= (t-m)\delta (A)$ and $X$
 is the direct sum of all indecomposable
 projective-injective $A^{(t)}$-modules which are not
 $A^{(m)}$-modules. Since ${\rm gl.dim} \ A^{(t)}\leq
 2t+1$, we know that $T'= T\oplus X$
 is a cotilting $A^{(t)}$-module. By [3], $^{\perp}T'= \{ M\in A^{(t)}-{\rm mod} \
 |\ {\rm Ext}^{i}_{A^{(t)}}(M,T)=0, i> 0\}$ is a contravariantly
 finite subcategory of $A^{(t)}$-mod.

Let $\mathcal {X}= {^{\perp}T}=\{ M\in A^{(m)}-{\rm mod} \
 |\ {\rm Ext}^{i}_{A^{(m)}}(M,T)=0, i> 0\}$. Then $\mathcal
 {X}= {^{\perp}T'}\cap A^{(m)}-{\rm mod}$.

We claim that  $\mathcal {X}$ is a contravariantly
 finite subcategory of $A^{(m)}$-mod.

In fact, $\forall\ M\in A^{(m)}-{\rm mod}$. Let $f_M: X_M\rightarrow
M$ be the minimal right $^{\perp}T'$-approximation of $M$. Then
$X_M\in A^{(m)}-{\rm mod}$ and $f_M: X_M\rightarrow M$ is the
minimal right $\mathcal {X}$-approximation of $M$. For otherwise,
$X_M$ can be decomposed as $X_M= X_m\oplus X_t$, where $X_m\in
A^{(m)}-{\rm mod}$ and $0\neq X_t\not\in A^{(m)}-{\rm mod}$. Since
$t> m$ and $X_t\in A^{(t)}-{\rm mod}$, we have that ${\rm
Hom}_{A^{(t)}}(X_t, M)=0$, which contradicts with the assumption
that $f_M$ is minimal.

By Lemma 2.3, we know that $T$ admits a complement $N$
 in mod-$A^{(m)}$, i.e., $T\oplus N$ is a tilting
$A^{(m)}$-module and $\delta (T\oplus N)= \delta
(A^{(m)})=(m+1)\delta (A)=\delta (T)$. This forces that ${\rm add} \
N\subset {\rm add} \ T$, therefore $T=T\oplus N$ is a tilting
$A^{(m)}$-module. This completes the proof.          \hfill$\Box$

 \vskip 0.2in

{\bf Theorem 3.2.} {\it Let $M$ be a partial tilting
$A^{(m)}$-module. Then $M$ admits a complement  $C$ in
$A^{(m)}$-mod.}

\vskip 0.1in

{\bf Proof.}  Let ${\rm pd}\ M = t$. If $t\leq m$, the consequence
has been proved in [2]. Therefore we only need to consider the case
of $m< t\leq 2m+1$.

By [2] again, we know that $M$ admits a complement $X$ in
$A^{(t)}$-mod with ${\rm pd}\ X \leq t$. Without loss generality, we
may assume that $M$ and $X$ are basic. We decompose $M$ as
$M=M_1\oplus P$, where $M_1$ has no projective-injective summand and
$P$ is a projective-injective $A^{(m)}$-module. We denote by $P_m$
the direct sum of all indecomposable projective-injective
$A^{(m)}$-modules and by $P_{(i)}$ the direct sum of all
indecomposable projective-injective $A^{(i)}$-modules which are not
$A^{(m)}$-modules for all $i$ with $m< i\leq t$. Then the tilting
$A^{(t)}$-module $M\oplus X$ can be written as $M\oplus X= M_1\oplus
P_m\oplus N_1\oplus P_{(t)}$ with $N_1$  haing no
projective-injective summand.

If ${\rm pd}\ N_1 \leq m$, then $N_1\in$ mod-$A^{(m)}$ and $\delta
(M_1\oplus P_m\oplus N_1) = (m+1)\delta (A)$. By Theorem 3.1, we
know that $M_1\oplus P_m\oplus N_1= M\oplus C$ is a tilting
$A^{(m)}$-module and that $C$ is a complement to $M$ in
$A^{(m)}$-mod.

If ${\rm pd}\ N_1 =i > m$, then $m< i\leq t$. Take an indecomposable
summand $N_j$ of $N_1$ with ${\rm pd}\ N_j =j > m$, let $N_{(j)}=
N_1/N_j$. Then $N_j$ is a complement to the faithful almost tilting
$A^{(t)}$-module $M_1\oplus P_m\oplus N_{(j)}\oplus P_{(t)}$.
According to Lemma 2.2 and Lemma 2.4, we have an exact sequence
$0\rightarrow X_j\rightarrow T_0\rightarrow \cdots \rightarrow
T_{j-1}\rightarrow N_j\rightarrow 0$ with ${\rm pd}\ X_j \leq m$ and
all $T_s\in {\rm add}(M_1\oplus P_m\oplus N_{(j)}\oplus P_{(t)})$
with $0\leq s\leq j-1$. Moreover, we know that $X_j$ is a complement
to the faithful almost tilting $A^{(t)}$-module $M_1\oplus P_m\oplus
N_{(j)}\oplus P_{(t)}$, and we obtain a tilting $A^{(t)}$-module
$M_1\oplus P_m\oplus N_{(j)}\oplus X_j\oplus P_{(t)}$.

We repeat the above process for every indecomposable summand with
projective dimension bigger than $m$. We then obtain a tilting
$A^{(m)}$-module $M\oplus P_m\oplus N'$, and we deduce that $M$ has
a tilting complement in $A^{(m)}$-mod. This completes the proof.
    \hfill$\Box$

\vskip 0.2in

{\bf Corollary 3.3.} {\it Let $A$ be representation-infinite and let
T be a faithful almost complete tilting $A^{(m)}$-module. Then the
number of non-isomorphic indecomposable complements to $T$ is either
$2m+1$ or $2m+2$.}

\vskip 0.1in

{\bf Proof.} Note that
$\Sigma_{2m}=\Omega_{A^{(m)}}^{-2m}\,\Sigma_0\subseteq \,A_m-{\rm
ind}$,  and $\Sigma_{2m-1}\cap A_m$-ind = $\emptyset$. If ${\rm pd}
T = 2m$, then the consequence can be deuced from Theorem 3.1 and
Lemma 2.1 and Lemma 2.2. Let ${\rm pd} T = 2m+1$. We denote by $P_m$
the direct sum of all indecomposable projective-injective
$A^{(m)}$-modules and by $P_{(i)}$ the direct sum of all
indecomposable projective-injective $A^{(i)}$-modules which are not
$A^{(m)}$-modules for all $i$ with $m< i\leq 2m+1$.  Then $T\oplus
P_{(2m+1)}$ can be regard as a faithful almost complete tilting
$A^{(2m+1)}$-module. By Lemma 2.2 and Theorem 3.1, we know that
$T\oplus P_{(2m+1)}$ has $2m+2$ non-isomorphic complements with
projective dimension at most $2m+1$, and that there are at least
$2m$ non-isomorphic complements with projective dimension at most
$2m$, these complements are $A^{(m)}$-modules.

Now if $T\oplus P_{(2m+1)}$ has only one complement with projective
dimension $2m+1$, then the number of non-isomorphic complements to
in $A^{(m)}$-mod to $T\oplus P_{(2m+1)}$ is at least $2m+1$, hence
$T$ has $2m+1$ non-isomorphic complements in $A^{(m)}$-mod.

If $T\oplus P_{(2m+1)}$ has exactly two complements $X_1$ and $X_2$
with ${\rm pd}_{A^{(2m+1)}} X_1={\rm pd}_{A^{(2m+1)}} X_2=2m+1$, the
by using Lemma 2.2, we have a non split exact sequence $0\rightarrow
X_1\stackrel{f}\rightarrow T'\stackrel{g}\rightarrow X_2\rightarrow
0$ with $T'\in {\rm add}\ T\oplus P_{(2m+1)}$ and $f:X_1\rightarrow
T'$ is the minimal left ${\rm add} \ T\oplus
P_{(2m+1)}$-approximation.  Since ${\rm pd}_{A^{(2m+1)}} X_1 = {\rm
pd}_{A^{(2m+1)}} X_2 = 2m+1$, $T'$ cannot be projective-injective.
It is easy to see that $0\neq T_0=T'/P_{2m+1}$ is a $A^{(m)}$-mod.
Since ${\rm Hom}_{A^{(2m+1)}}(X_1, T_0)\neq 0$, $T_0$ must have an
indecomposable summand $X$ such that ${\rm Hom}_{A^{(2m+1)}}(X_1,
X)\neq 0$. By Lemma 2.1, we have that ${\rm Hom}_{A^{(m)}}(X_1,
X)\neq 0$, and that $X_1$ is a predecessor of $X$. Note that $X$ is
a $A^{(m)}$-module,  therefore $X_1\in A^{(m)}$-mod, and $T$ has
$2m+1$ non-isomorphic complements in $A^{(m)}$-mod. This completes
the proof. \hfill$\Box$

\vskip 0.2in

Let $\mathcal {T}_{A^{(m)}}$ be the set of all basic tilting
$A^{(m)}$-modules up to isomorphism. According to [10], the tilting
quiver $\overrightarrow{\mathscr{K}}_{A^{(m)}}$ of tilting
$A^{(m)}$-modules is defined as the following. The vertices of
$\mathscr{K}_{A^{(m)}}$ are the elements of $\mathcal
{T}_{A^{(m)}}$. There is an arrow $T'\rightarrow T$ if
$T'=\overline{T}\oplus X$, $T=\overline{T}\oplus Y$ with $X, Y$
indecomposable and there is a short exact sequence $0\rightarrow
X\rightarrow E\rightarrow Y\rightarrow 0$ with $E\in{\rm add}\
\overline{T}$. The underlying graph of
$\overrightarrow{\mathscr{K}}_{A^{(m)}}$ is denoted by
$\mathscr{K}_{A^{(m)}}$, and it is called the tilting graph of
$A^{(m)}$.

\vskip 0.2in

 {\bf Theorem 3.4.} {\it Let $A^{(m)}$ be the $m$-replicated algebra
 of a hereditary algebra of $A$. Then the tilting quiver
 $\overrightarrow{\mathscr{K}}_{A^{(m)}}$ of  $A^{(m)}$ is connected.}

\vskip 0.1in

{\bf Proof.} Let $T$ be a basic tilting $A^{(m)}$-module with ${\rm
pd} \ T = t\geq 1$, and $X_t$ be an indecomposable summand of $T$
with ${\rm pd} \ X_t = t$. Then $T_{(t)}= T/ X_t$ is a faithful
almost tilting $A^{(m)}$-module and $X_t$ is a complement to
$T_{(t)}$. By Lemma 2.2, we have the following exact sequence.
$$ (*)\ \ \ \ \   0\longrightarrow
X_0\stackrel{f_{0}}\longrightarrow
T_{0}^\prime\stackrel{f_{1}}\longrightarrow T_{1}^\prime
\longrightarrow\cdot\cdot\cdot\longrightarrow
T_{t-2}^\prime\stackrel{f_{t-1}}\longrightarrow
T_{t-1}^\prime\stackrel{f_{t}}\longrightarrow X_t\longrightarrow 0$$
in $A^{(m)}$-mod,  such that

\vskip 0.1in

{(1)} $T_i^\prime\in{\rm add}\, T_{(t)} $ for all $0\leq i\leq t-1$,

\vskip 0.1in

{(2)} $ X_i={\rm Im}\ f_i$ for  $0\leq i\leq t$,  and $i\leq {\rm
pd}_{A^{(m)}}X_i\leq i+1$
 for  $0\leq i\leq t-1$,

\vskip 0.1in

{(3)} each of the induced monomorphisms $X_i\hookrightarrow
T_i^\prime$ is a  minimal left  ${\rm add}\ T_{(t)}$-approximation.

Note that ${\rm pd}_{A^{(m)}}X_0\leq 1$, and by $(*)$ there is a
path from $X_0\oplus T_{(t)}$ to $T=T_{(t)}\oplus X_t$ in
$\overrightarrow{\mathscr{K}}_{A^{(m)}}$. We can repeat this process
for indecomposable summands with maximal projective dimension until
obtaining a path from a tilting $A^{(m)}$-module with projective
dimension  at most one to $X_0\oplus T_{(t)}$. According to [15] we
know that all basic tilting  $A^{(m)}$-modules with projective
dimension  at most one consists of a connected subquiver of
$\overrightarrow{\mathscr{K}}_{A^{(m)}}$, so
$\overrightarrow{\mathscr{K}}_{A^{(m)}}$ is connected. This
completes the proof.       \hfill$\Box$

\vskip 0.2in

\section {Complements to almost complete tilting modules over duplicated algebras}

\vskip 0.2in

In this section, we investigate the number of complements to a basic
almost complete tilting module over duplicated algebras. According
to [6], we know that if a basic almost complete tilting module is
not faithful, then it has a unique complement. Therefore we only
need to consider faithful basic almost complete tilting modules over
duplicated algebras.

From now on,  we always assume that $A$ is a representation-infinite
hereditary algebra over an algebraically closed field $k$, and
$A^{(1)}$ is the duplicated algebra of $A$. Note that ${\rm gl.dim}\
A^{(1)}=3$. According to Corollary 3.3, we know that a faithful
basic almost complete tilting modules over duplicated algebras has
at least three complements and at most four complements.

\vskip0.2in

{\bf Proposition 4.1.} {\it Let $A^{(1)}$ be the duplicated algebra
 of $A$ and $T_1$ be a faithful basic almost
 complete tilting $A^{(1)}$-module with ${\rm pd}_{A^{(1)}}T_1\leq 1$.
 If $T_1$ has exactly four non-isomorphic complements, then
 $T_0= T_1/P$ is a faithful $A_0$-module, where $P_1$ is the direct
 sum of projective injective summands of $T_1$.}

\vskip 0.1in

{\bf Proof.} \ Assume that $T_1$ has four non-isomorphic complements
$X_0, X_1, X_2, X_3$. Since ${\rm pd}_{A^{(1)}}T_1\leq 1$, by Lemma
2.2 and [15], we can assume that ${\rm pd}\ X_0\leq 1$ and ${\rm
pd}\ X_1\leq 1$, ${\rm pd}\ X_2 = 2$  and ${\rm pd}\ X_3 = 3$. If
$T_0= T_1/P_1$ is not a faithful $A_0$-module, according to Lemma
2.1 we can regard $A^{(1)}$-ind as a full convex subcategory of
$A^{(3)}$-ind. By [7] and [15], we know that $X_1\in \Sigma_1$,
$X_2\in \Sigma_2$ and $X_3\in \Sigma_3$. Since $\Sigma_{3}\cap
A_1$-ind = $\emptyset$, we have that $X_3\not\in A^{(1)}$-mod. This
contradiction insures that $T_0= T_1/P_1$ must be a faithful
$A_0$-module. The proof is completed. \hfill$\Box$

\vskip 0.2in

{\bf Remark.} \ Generally speaking, the converse of Proposition 4.1
dose not holds.

\vskip 0.2in

{\bf Example 1.}\  Let $\widetilde{D}_4$ be the tame quiver and
$A^{(1)}=k\widetilde{D_4}^{(1)}/I$.

$$\begin{array}{ccccccc}

       && 2 & & &&2'\\[-2ex]

       &\swarrow & &\nwarrow & &\swarrow& \\[-2ex]

  \widetilde{D_4}^{(1)}:    1 &\leftleftarrows & \begin{array}{c}3\\[-2ex]4\end{array}
    &\leftleftarrows & 1' &\leftleftarrows &  \begin{array}{c}3'\\[-2ex]4'\end{array}\\[-2ex]

      &  \nwarrow & & \swarrow&&\searrow &\\[-2ex]

    && 5 &&&&5'\\[-2ex]

    \end{array}.
$$

\vskip 0.2in

Then the indecomposable projective-injective $A^{(1)}$-modules are
represented by their Loewy series as the following,
$$
P_1'=\begin{array}{c}1'\\[-2ex]2345\\[-2ex]1\end{array},\
P_2'=\begin{array}{c}2'\\[-2ex]1'\\[-2ex]2\end{array},\
P_3'=\begin{array}{c}3'\\[-2ex]1'\\[-2ex]3\end{array}, \
P_4'=\begin{array}{c}4'\\[-2ex]1'\\[-2ex]4\end{array}, \
P_5'=\begin{array}{c}5'\\[-2ex]1'\\[-2ex]5\end{array}.
$$

The direct sum $P'_1\oplus P'_2\oplus P_3'\oplus P'_4\oplus P_5'$ of
all indecomposable projective-injective $A^{(1)}$-modules is denoted
by $P_1$.

\vskip 0.1in

 Let $T_1=\begin{array}{c}2345\\[-2ex]1
\end{array}\oplus 3\oplus 4\oplus 5\oplus P_1$. Then
 $T_1$ is a faithful basic almost complete tilting $A^{(1)}$-module with
 ${\rm pd}\ T_1 =1$ and $T_0=T_1/P_1$ is a faithful
$A_0$-module.

One can easily see that $T_1$ has only three non-isomorphic
complements $X_0=\begin{array}{c}345\\[-2ex]1\end{array}$,
$X_1=2$, $X_2= \begin{array}{c}2'\\[-2ex]1'\end{array}$,
and ${\rm pd}\ X_0 ={\rm pd}\ X_1=1$, ${\rm pd}\ X_2=2$.  Note that
$T_1$ has no complement with projective dimension 3.

\vskip0.2in

{\bf Theorem 4.2.} {\it Let $A^{(1)}$ be the duplicated algebra
 of $A$, and $T_1$ be a faithful basic almost
 complete tilting $A^{(1)}$-module with ${\rm pd}_{A^{(1)}}T_1\leq 1$.
 Then $T_1$ has exactly four non-isomorphic complements if and only if
 the unique complement $X$ to $T_1$ with ${\rm pd}_{A^{(1)}} X =2$ such that the
 injective envelope $E(X)$ of $X$ is projective.}

\vskip 0.1in

{\bf Proof.} \ We assume that the unique complement $X$ to $T_1$
with ${\rm pd}_{A^{(1)}} X =2$ such that the injective envelope
$E(X)$ of $X$ is projective, then we have an exact sequence
$0\rightarrow X \stackrel{f}\rightarrow E(X)\stackrel{g}\rightarrow
Y\rightarrow 0$ with $Y\neq 0$. By Lemma 6 in [2],  we know that
$E(X)$ is the projective cover of $Y$ and $Y$ is also
indecomposable. Note that ${\rm Ext}^1_{A^{(1)}}(T_1, X)=0$, we
deduce that $g: E(X)\rightarrow Y$ is the right minimal ${\rm add} \
T_1$-approximation. By using Lemma 2.2, we know that $Y$ is also a
complement to $T_1$ with ${\rm pd}\ Y =3$, hence $T_1$ has four
non-isomorphic complements.

\vskip 0.1in

Conversely,   if $T_1$ has four non-isomorphic complements, then
$T_1$ has two complements with projective dimension at most one, a
unique complement $X$ with ${\rm pd}\ X=2$ and a unique complement
$Y$ with ${\rm pd}\ Y =3$ respectively. The projective cover $P(Y)$
of $Y$ is injective since $Y\in A^{(1)}$-mod. According to Lemma 6
in [2], we have a non split exact sequence $0\rightarrow K
\stackrel{f}\rightarrow P(Y)\stackrel{g}\rightarrow Y\rightarrow 0$
with $K$ indecomposable and ${\rm pd}\ K =2$. Moreover, $P(Y)$ is
the injective envelope of $K$. Note that $g: P(Y)\rightarrow Y$ is
the right minimal ${\rm add} \ T_1$-approximation since ${\rm
Ext}^1_{A^{(1)}}(T_1, K)=0$. By using Lemma 2.2, we know that $K$ is
a complement to $T_1$ with ${\rm pd}\ K =2$. It forces that $X\simeq
K$ and the injective envelope $E(X)= P(Y)$ of $X$ is projective.
This completes the proof.        \hfill$\Box$

\vskip 0.2in

{\bf Theorem 4.3.} {\it Let $A^{(1)}$ be the duplicated algebra
 of $A$ and $T_2$ be a faithful basic almost
 complete tilting $A^{(1)}$-module with ${\rm pd}_{A^{(1)}}T_2 \leq 2$.
 Then $T_2$ has exactly four non-isomorphic complements if and only if
  $T_2$ has a complement $X$ with ${\rm pd}\ X =3$.}

\vskip 0.1in

{\bf Proof.} \ If $T_2$ has a complement $X$ with ${\rm pd}\ X =3 $,
 then by Lemma 2.2, $T_2$ has three non-isomorphic complements with
 projective dimension at most 2, hence $T_2$ has exactly four non-isomorphic complements.

 Conversely, if $T_2$ has exactly four non-isomorphic complements,
 by using Lemma 2.2 again, we know that $T_2$ must have a
 complement $X$ with ${\rm pd}\ X =3$.       \hfill$\Box$

\vskip 0.2in

{\bf Corollary 4.4.} {\it Let $A^{(1)}$ be the duplicated algebra
 of $A$ and $T_2$ be a faithful basic almost complete tilting $A^{(1)}$-module
 with ${\rm pd}_{A^{(1)}}T_2 =2$. If $T_2$ has a complement $X$
 with ${\rm pd}\ X =2$ such that the injective envelope $E(X)$ of $X$ is
 projective, then $T_2$ has exactly four non-isomorphic complements.}

\vskip 0.1in

{\bf Proof.}\ According to the assumption, we have a non split exact
sequence $0\rightarrow X \stackrel{f}\rightarrow
E(X)\stackrel{g}\rightarrow Y\rightarrow 0$ with $Y$ indecomposable
and ${\rm pd}\ Y =3$ and $E(X)$ is the projective cover of $Y$.
Since ${\rm Ext}^1_{A^{(1)}}(T_2, X)=0$, we deduce that $g:
E(X)\rightarrow Y$ is the right minimal ${\rm add} \
T_2$-approximation. By using lemma 2.2, we know that $Y$ is a
complement to $T_2$. Since ${\rm pd}\ Y =3$, by Proposition 4.3 we
know that $T_2$ has exactly four non-isomorphic complements. This
completes the proof. \hfill$\Box$

\vskip 0.2in

{\bf Remark 4.5.} \ Note that the converse is not true. That is, let
$T_2$ be a faithful basic almost complete tilting $A^{(1)}$-module
with ${\rm pd}_{A^{(1)}}T_2 =2$. If $T_2$ has four non-isomorphic
complements, we know that $T_2$ must have complements with
projective dimension 2, but their injective envelope may be not
projective.

\vskip 0.2in

{\bf Example 2.}\  Let $A^{(1)}=k\widetilde{D_4}^{(1)}/I$ be the
same as example 1.

Let $T_2=\begin{array}{c}3'4'5'\\[-2ex]1'\ 1'
\end{array}\oplus \begin{array}{c}2'4'5'\\[-2ex]1'\ 1' \end{array}\oplus
\begin{array}{c}2'3'5' \\[-2ex]1'\ 1' \end{array} \oplus
\begin{array}{c}2'3'4' \\[-2ex]1'\ 1'\end{array}\oplus P_1$.

Then
 $T_2$ is a faithful basic almost complete tilting $A^{(1)}$-module with
 ${\rm pd}\ T_2 =2$. One can easily see that $T_2$ has four
nonisomorphic complements $X_0=1$, $X_1=
\begin{array}{c}1'\\[-2ex]2345\end{array}$, $X_2= \begin{array}{c}2'3'4'5'\\[-2ex]1'1' 1'
\end{array}$, $X_3= \begin{array}{c}2'2'3'3'4'4'5'5'\\[-2ex]1' 1' 1' 1' 1'
\end{array}$, and ${\rm pd}\ X_i =i$ for $0\leq i\leq 3$. However the injective
envelope $E(X_2)= (\begin{array}{c}2'3'4'5'\\[-2ex]1'
\end{array})^3$ of $X_2$ is not projective.

\vskip 0.2in

{\bf Remark 4.6.} \ We should mention that Theorem 4.3 dose not hold
for a faithful basic almost complete tilting $A^{(1)}$-module $T_3$
with ${\rm pd}_{A^{(1)}}T_3 =3$. That is,  if $T_3$ has  a
complement $X$ with ${\rm pd}_{A^{(1)}}X =3$, we cannot deduce that
$T_3$ has four non-isomorphic complements.

\vskip 0.2in

{\bf Example 3.}\  Let $Q: 1\leftleftarrows 2\leftleftarrows
1'\leftleftarrows 2'$ and $A^{(1)}= kQ/I$ be the duplicated algebra
of the Kronecker algebra. The indecomposable projective-injective
$A^{(1)}$-modules are
$P_1'= \begin{array}{c}1'\\[-2ex]22\\[-2ex]1\end{array}$ and
$P_2'= \begin{array}{c}2'\\[-2ex]11\\[-2ex]2\end{array}$
which are represented by their Loewy series.

\vskip 0.1in

(1). Let $T_3= 2'\oplus  \begin{array}{c}1'\\[-2ex]22\\[-2ex]1\end{array}
\oplus \begin{array}{c}2'\\[-2ex]11\\[-2ex]2\end{array}$. Then
${\rm pd}_{A^{(1)}}T_3 =3$, and $T_3$ has only three non-isomorphic
complements $X_1=\begin{array}{c}2222\\[-2ex]111 \end{array}$,
$X_2=\begin{array}{c}1'1'1'\\[-2ex]22 \end{array}$,
$X_3=\begin{array}{c}2'2'\\[-2ex]1' \end{array}$, and ${\rm pd}_{A^{(1)}}X_i
=i$ for $1\leq i\leq 3$.

\vskip 0.1in

(2). The following examples indecate that a faithful basic almost
complete tilting $A^{(1)}$-module $T_2$ with ${\rm pd}_{A^{(1)}}T_2
\leq 2$ may have no complement with projective dimension 3.

\vskip 0.1in

(i)\ Let $T_1= 2\oplus \begin{array}{c}1'\\[-2ex]22\\[-2ex]1\end{array}
\oplus \begin{array}{c}2'\\[-2ex]11\\[-2ex]2\end{array}$. Then
${\rm pd}_{A^{(1)}}T_1 =1$, and $T_1$ has only three non-isomorphic
complements $X_0=\begin{array}{c}22\\[-2ex]1 \end{array}$,
$X_1=\begin{array}{c}1'\\[-2ex]22 \end{array}$,
$X_2=1'$, and ${\rm pd}_{A^{(1)}}X_0 ={\rm pd}_{A^{(1)}}X_1=1$,
${\rm pd}_{A^{(1)}}X_2 =2$. Note that $T_1$ has no complement with
projective dimension 3.

\vskip 0.1in

(ii)\ Let $T_2= \begin{array}{c}2'\\[-2ex]1'1'\\[-2ex]1\end{array}
\oplus \begin{array}{c}1'\\[-2ex]22\\[-2ex]1\end{array}
\oplus \begin{array}{c}2'\\[-2ex]11\\[-2ex]2\end{array}$. Then
${\rm pd}_{A^{(1)}}T_2 =2$, and $T_2$ has  three non-isomorphic
complements $X_1=2$, $X_2=1'$, $X_2'=\begin{array}{c}2'2'\\[-2ex]1'1'1' \end{array}$,
 and ${\rm pd}_{A^{(1)}}X_1 =1$,  ${\rm pd}_{A^{(1)}}X_2={\rm pd}_{A^{(1)}}X_2' =2$.
 Note that $T_2$ has no complement with projective dimension 3.

\vskip 0.2in

Let $\overline{T}$ be a faithful basic almost
 complete tilting $A^{(1)}$-module. If ${\rm pd}\ \overline{T} \leq 1$, then
$\overline{T}$ has a complement with projective dimension 1, and if
${\rm pd}\ \overline{T} = 2$,  $\overline{T}$ always has a
complement with projective dimension 2.  Now let $T_3$ be a faithful
basic almost
 complete tilting $A^{(1)}$-module with ${\rm pd}_{A^{(1)}}T_3 =3$.
 We know that $T_3$  has complements $X_1$
 and $X_2$ in $A^{(1)}$-{\rm mod} with ${\rm pd}\ X_1 \leq 1$
 and ${\rm pd}\ X_2 =2$ respectively. An interesting question is
 that whether $T_3$ also has a complement $X_3$
 in $A^{(1)}$-{\rm mod} with  ${\rm pd}\ X_3 =3$.

\begin{description}

\item{[1]}\ I.Assem, T.Br$\ddot{{\rm u}}$stle, R.Schiffler, G.Todorov,
   Cluster categories and duplicated algebras.
   {\it J. Algebra }, 305(2006), 548-561.

\item{[2]}\ I.Assem, T.Br$\ddot{\rm u}$stle, R.Schiffer, G.Todorov,
 $m$-cluster categories and $m$-replicated algebras.  {\it Journal
of pure and applied algebra},  212(2008), 884-901.

\item{[3]}\ M.Auslander, I.Reiten,  Applications of contravariantly finite
subcategories.  {\it Adv. Math.}, 86(1991), 111-152.

\item{[4]}\ M.Auslander, I.Reiten, S.O.Smal$\phi$,  Representation
Theory of Artin Algebras.  Cambridge Univ. Press, 1995.

\item{[5]}\ K. Bongartz, {\it Lecture Notes in Math. 903.}
Springer-Verlag, Berlin, Heidelberg, New York, 1981.

\item{[6]}\ F.Coelho, D.Happel, L.Unger,   Complements to partial tilting
modules. {\it J. Algebra }, 170(1994), 184-205.

\item{[7]}\  D.Happel, L.Unger, Almost complete tilting modules.
 {\it Proc.Amer.Math. Soc.}, 107(1989), 603-610.

\item{[8]}\  D.Happel, L.Unger,  Partial tilting modules and covariantly finite
subcategories. {\it Comm.Algebra},  22(1994), 1723-1727.

\item{[9]}\ D.Happel, L.Unger,  Complements and the generalized Nakayama
conjecture.
{\it CMS Conf. Proc.} 24(1998), 293-310.

\item{[10]}\ D.Happel, L.Unger, On the quiver of tilting modules. J.Algebra, 284(2005), 857-868.

\item{[11]}\ X.Lei, H.Lv, S.Zhang, Complements to the almost complete tilting
$A^{(m)}$-modules.  {\it Comm.Algebra},  37(2009), 1719-1728.

\item{[12]}\ H.Lv, S.Zhang, Representation dimension of $m$-replicated
algebras. Preprint, arXiv:math.Rt/0810.5185, 2008.

\item{[13]}\ J.Rickard, A. Schofield,  Cocovers and tilting modules.
{\it Math. Proc. Cambridge Philos. Soc.},  106(1985), 1-5.

\item{[14]}\ C.M.Ringel,  Tame algebras and integral quadratic forms.
{\it Lecture Notes in Math. 1099.}  Springer-Verlag, Berlin,
Heidelberg, New York, 1984.

\item{[15]}\ S.Zhang,  Tilting mutation and duplicated
algebras.  {\it Comm.Algebra}, 37(2009), 3516-3524.

\end{description}

\end{document}